\documentclass[11pt]{article}

\usepackage{mathrsfs}
\usepackage{bbm}
\usepackage{amsthm}
\usepackage{amsmath,amssymb}

\usepackage{graphicx}
\usepackage{cases}

\newcommand{\head}[0]{\mathfrak{h}}
\newcommand{\tail}[0]{\mathfrak{t}}
\newcommand{\inv}[0]{^{-1} }
\newcommand{\arc}[0]{\mathcal{A}}

\newcommand{\mat}[0]{\mbox{\rm Mat}}

\title {A weighted graph zeta function \\ involved in the Szegedy walk}
\author{Ayaka Ishikawa \\
Graduate School of Engineering Science\\
Yokohama National University \\ 
Hodogaya, Yokohama, 240-8501, JAPAN \\
Norio Konno \\
Department of Applied Mathematics, Faculty of Engineering \\ 
Yokohama National University \\
Hodogaya, Yokohama, 240-8501, JAPAN}
\date{}
\begin{document}
\maketitle

\begin{abstract}
We define a new weighted zeta function for a finite graph
and obtain its determinant expression.
This result gives the characteristic polynomial of the transition matrix of the Szegedy walk on a graph.

\end{abstract}


\section{Introduction}

In this paper, we define a new graph zeta function which is not obtained as a specialization of the generalized weighted zeta function.
Since it satisfies the adjacency condition, we can obtain the three expressions.
The problem to consider in this paper is to give the Ihara expression.
In the main theorem, we show that our graph zeta function has the Ihara expression for a finite graph and a finite digraph.
The reason we deal with graphs is to consider the relationship between our graph zeta function and the ``Szegedy walk'' defined on a graph.
On the other hand, we also treat digraphs.
Since the Ihara expressions of generalized weighted zeta functions on graphs and digraphs are a slightly different (cf. \cite{IDE2021227, ishikawa21}),
we can expect there are differences between our graph zeta functions for graphs and digraphs.
In fact, the two Ihara expressions are find to be different.
The reason lies in the difference in the definition of the inverse arc as expected.
For a graph,
we consider the symmetric digraph of the graph
and define a pair of arcs assigned to each edge as inverse.
Then the Ihara expression is the same form as in \cite{IDE2021227}.
We deal with this definition because the Szegedy walk adopts the definition. 
In other words, it is necessary to use the definition
in order to consider the relationship between the Szegedy walk and our graph zeta function.
In this paper, we deal with the other definition of the inverse motivated by \cite{ishikawa21}.
For a digraph and an arc,
we define the inverse as any arc with an opposite direction.
Then the Ihara expression different from that derived by the previous definition.

The Szegedy walk is a quantum walk model, which is a generalization of the ``Grover walk''.
A quantum walk is the quantum version of a random walk,
and it is studied in various fields: 
quantum algorithm, financial engineering, and laser isotope separation, for example (see e.g., \cite{matsuoka2011theoretical, Orrell_2020, Portugal_2018}).
The common interest in these fields is the behavior of quantum walks.
Some important quantum walk properties, 
such as periodicity and localization, are determined by the eigenvalues of the transition matrix.
For example,
in the quantum search problem,
the eigenvalues are crucial for deriving the quantum hitting time of the quantum walk on a graph \cite{konno2017}.
In particular, 
the Grover transition matrix is an example of  the ``edge matrix'' \cite{morita2020ruelle}, which describes the Hashimoto expression for the Sato zeta function.
By transforming the Hashimoto expression of the Sato zeta function into the Ihara expression,
this process describes the eigenvalues of the Grover transition matrix \cite{konno2012relation}.
Our graph zeta function has the edge matrix corresponding to the Szegedy transition matrix.
Hence, our result shows that we can describe the eigenvalues of the Szegedy transition matrix by our graph zeta function.

The rest of the paper is organized as follows.
In Section \ref{sec:pre}, we define a new graph zeta function as an exponential expression.
We introduce the Szegedy walk and state the relation between our zeta function and the Szegedy walk.
The Ihara expression of our zeta function is defined in Section \ref{main} which is the main theorem of our paper.
From the theorem,
we can obtain the characteristic polynomial of the transition matrix of the Szegedy walk on a finite graph.
In Section \ref{sec:exam}, we give two examples on a graph and a digraph.

Throughout this paper, 
graphs (resp. digraphs) are finite, and multi-edge (resp. multi-arcs) and multi-loops are allowed.
We use the following symbols.
For positive integers $m$ and $n$, let $\mat(m,n;\mathbb{C})$ be the set of $m\times n$ matrices over $\mathbb{C}$.
For a matrix $M\in\mat(n,n;\mathbb{C})$,
let ${\rm Spec}(M)$ be the spectrum of $M$.
An $m\times n$ matrix with all one denotes by $\mathbbm{1}_{m\times n}$.
In particular, if $m=n$, we write $\mathbbm{1}_{m\times n}$ as $\mathbbm{1}_{n}$.
For a proposition $P$,
we define $\delta_{P}$ as follows: $\delta_{P}=1$ if $P$ is true, $\delta_{P}=0$ if $P$ is false.
For the Kronecker delta $\delta_{uv}$,
let $\overline{\delta}_{uv}=1$ if $u\neq v$, $\overline{\delta}_{uv}=0$ otherwise.

\section{Preliminary}\label{sec:pre}
\noindent
A {\it graph} $G=(V,E)$ is a pair of {\it vertex} set $V$ and {\it edge} set $E$, 
where $E$ consists of $2$-subsets of $V$.
If both $V$ and $E$ are finite, then $G$ is called {\it finite}.
We call an edge $e=\{v,v\}$ a loop.
The number $\deg(u) := \#\{ \{u,v\} \in E | v\in V \}$ is called the {\it degree} of $u$.
If there is at most one edge between each two vertices and there are not loops,
then the graph is called {\it simple}.
Let $V$ be a vertex set and $\arc$ a set of ordered pairs of two vertices.
We call the pair $\Delta =(V,\arc)$ a {\it digraph} and an element of $\arc$ an {\it arc}.
For an arc $a=(u,v)$, 
$u$ and $v$ are called the {\it tail} and the {\it head} of $a$ denoted by $\tail(a)$ and $\head(a)$, respectively.
For two vertices $u,v\in V$ of a digraph $\Delta$, 
let $\arc_{uv}:= \{ a\in\arc \ | \ \tail(a)=u, \head(a)=v \}$, $\arc_{u*} := \{ a\in\arc | \tail(a)=u\}$, $\arc_{*v} := \{ a\in\arc \ | \ \head(a)=v\}$,
and $\arc(u,v) := \arc_{uv}\cup \arc_{vu}$.
For a graph $G=(V,E)$,
let $\arc(G):=\{ (u,v), (v,u) | e=\{ v,u \}\in E \}$,
and then the digraph $\Delta(G)=(V,\arc(G))$ is called the {\it symmetric digraph} of $G$.
For an arc $a\in\arc(G)$, we denote by $\overline{a}$ the arc induced by the same edge as $a$.

A sequence of arcs $p=(a_i)_{i=1}^{k}$ is a {\it path} 
if it satisfies $\head(a_i) = \tail(a_{i+1})$ for each $i=1,2,\ldots, k-1$.
The number $k$, called the {\it length} of $p$, is denote by $|p|$.
If $\head(a_k)=\tail(a_1)$, then the path $p$ is called {\it closed}.
Let $X_k$ denote the set of closed paths of length $k$.
For $C\in X_k$,
we denote by $C^n$ the closed path that connects $C$ $n$ times.
It is called the $n$-th power of $C$.
If $C$ cannot be expressed as a power of a closed path shorter than $C$, then it is called {\it prime}.
For $C=(c_i)_{i=1}^k, C'=(c_i')_{i=1}^k \in X_k$, 
if there exists an integer $n$ such that $c_i=c'_{i+n}$ for any $i$, where the indices are taken modulo $k$,
then we denote the relation by $C\sim C'$.
Clearly, the relation $\sim$ is an equivalence relation.
An equivalence class is called a {\it cycle}, and we denote by $[C]$ the equivalence class of a closed path $C$.
Since any closed path in $[C]$ have the same length, 
we define the {\it length} of $[C]$ to be the length of a closed path in $[C]$.
We denote by $|C|$ the length of $[C]$.
A cycle is {\it prime} if a closed path in the cycle is prime.
We denote by $\mathcal{P}$ the set of prime cycles.


\subsection{A new graph zeta function}

Let $\Delta=(V,\arc)$ be a digraph.
For a map $\theta : \arc\times\arc \to \mathbb{C}$ and a closed path $C=(c_i)_{i=1}^k\in X_k$, 
let ${\rm circ}_\theta (C)$ denote the circular product 
$
	\theta(c_1,c_2) \theta(c_2,c_3) \ldots \theta(c_k,c_1).
$
Note that ${\rm circ}_\theta(C)={\rm circ}_\theta(C')$ holds if $C\sim C'$.
Let $N_k({\rm circ}_\theta) := \sum_{C\in X_k} {\rm circ}_\theta(C)$.
We define a graph zeta function for $\Delta$.

\vspace*{12pt}
\noindent
{\bf Definition~1:}

A graph zeta function for $\Delta$ is the following formal power series:
	\begin{align}\label{ident:zeta}
		Z_\Delta(t;\theta) := \exp \left(\sum_{k\ge 1}\frac{N_k({\rm circ}_\theta)}{k}t^k \right).
	\end{align}

\vspace*{12pt}
\noindent

We call the map $\theta$ the {\it weight} of the graph zeta function,
and the expression (\ref{ident:zeta}) the {\it exponential expression} \cite{morita2020ruelle}.
Let 
\begin{align*}
	E_\Delta (t;\theta) := \prod_{[C]\in\mathcal{P}} \frac{1}{1-{\rm circ}_\theta (C)t^{|C|}}, \qquad 
	H_\Delta (t;\theta) := \frac{1}{\det(I-tM_\theta)},
\end{align*}
where $M_\theta=(\theta(a,a'))_{a,a'\in\arc}$.
The expressions $E_\Delta(t;\theta)$ and $H_\Delta(t;\theta)$ are called
the {\it Euler expression} and the {\it Hashimoto expression},
respectively (cf. \cite{morita2020ruelle}).

\vspace*{12pt}
\noindent
{\bf Proposition~1:}
	
	If $\theta : \arc\times\arc \to \mathbb{C}$ satisfies the condition 
	$$
		\theta(a,a') \neq 0 \Rightarrow \head(a)=\tail(a'),
	$$
	then $Z_\Delta(t;\theta)= E_\Delta (t;\theta)=H_\Delta (t;\theta)$.

\vspace*{12pt}
\noindent
{\bf Proof:}
	See \cite{morita2020ruelle}.  $\square\,$

\vspace*{12pt}
\noindent

The above condition for $\theta$ is called the {\it adjacency condition} \cite{morita2020ruelle}.

Before we introduce the ``Ihara expression'',
we will mention the definition of the inverse.
For an arc $a$, let $a\inv$ denote the set of inverses of $a$.
For a digraph $\Delta$, if $a\in\arc_{uv}$, then we define $a\inv :=\arc_{vu}$.
Note that for a loop $a\in\arc_{vv}$, $a\in a\inv$ holds.
Unless otherwise specified, for a digraph, this definition is adopted.
On the other hand, for the symmetric digraph $\Delta(G)$ of a graph $G$,
we adopt the following definition.
Let $a=(u,v)$ and $a'=(v,u)$ be induced by an edge $\{ u,v\}\in E$,
then we define $a\inv :=\{ a' \}$, and $a'$ denotes by $\overline{a}$.

Let $\Delta=(V,\arc)$ be a digraph (not specified whether it is the symmetric digraph or not).
For any maps $\tau_1, \tau_2 : \arc \to \mathbb{C}$,
let $\tau$ be a map $\arc\times \arc \to \mathbb{C}$ defined by $\tau(a,a') := \tau_1(a)\tau_2(a')$.
We define a new graph zeta function $Z_\Delta(t;\theta)$ with the weight
\begin{align} \label{def:weight}
	\theta(a,a') := \tau(a,a')\delta_{\head(a)\tail(a')} - \delta_{a'\in a\inv}, 
\end{align}
where 
$\delta_{\head(a)\tail(a')}$ is the Kronecker delta, and
$\delta_{a'\in a\inv}$ equals $1$ if $a'$ is an inverse of $a$, $0$ otherwise.
Note that if $\delta_{a'\in a\inv}=1$,
any arc $a' \in a\inv$ is in $\arc_{\head(a) \tail(a)}$, and it satisfies $\head(a)=\tail(a')=v$.
If $\theta(a,a')\neq 0$ for $a,a'\in\arc$, then $\delta_{\head(a) \tail(a')} =1$ holds at least.
Thus, the weight $\theta$ is satisfying the adjacency condition,
and 
we can see $Z_\Delta(t;\theta)= E_\Delta (t;\theta)=H_\Delta (t;\theta)$.

We assume that $\theta^{\rm S} := \theta|_{\tau_1 \equiv 1}$.
The graph zeta function with the weight $\theta^{\rm S}$ is called the {\it Sato zeta function} \cite{morita2020ruelle, sato2007new}.
Since $\theta^{\rm S}$ satisfies the adjacency condition,
the Sato zeta function has the Euler expression and the Hashimoto expression.
In addition, it also has the following determinant expression called the {\it Ihara expression}.

\vspace*{12pt}
\noindent
{\bf Proposition~2:} (Sato \cite{sato2007new}) 
	Let $\Delta(G)=(V,\arc(G))$ be the symmetric digraph of a finite simple graph $G=(V,E)$,
	and $A^{\rm S}_G = (a_{uv})_{u,v\in V}$ and $D^{\rm S}_G =(d_{uv})_{u,v\in V}$
	be defined as follows:
	\begin{align*}
		a_{uv} := \sum_{a\in \arc_{uv}}\tau(a), \  \ \
		d_{uv} := \delta_{uv}\sum_{a\in \arc_{v*}} \tau(a).
	\end{align*}
	Then, $Z_{\Delta(G)}(t;\theta)$ equals
	$$
		\frac{1}{(1-t^2)^{|E|-|V|}\det(I-tA^{\rm S}_G +t^2 (D^{\rm S}_G -I))}.
	$$

	Let $\Delta=(V,\arc)$ be a digraph.
	We fix a total order $<$ on $V$.
	Let $\Phi_\Delta = \{ (u,v)\in V\times V \ | \ \arc(u,v)\neq \emptyset \}$.
	If $(u,v)\in\Phi_\Delta$, we assume $u<v$.
	Note that $\arc(u,v)=\arc(v,u)$ holds.
	Let $\underline{A}^{\rm S}_\Delta= (\underline{a}_{uv})_{u,v\in V}$ and $\underline{D}^{\rm S}_\Delta=(\underline{d}_{uv})_{u,v\in V}$
	be defined as follows:
	\begin{align*}
		\underline{a}_{uv} := \sum_{a\in \arc_{uv}} \frac{\tau(a)}{f(u,v)}, \  \ \
		\underline{d}_{uv} := \delta_{uv}\sum_{a\in \arc_{u*}} \frac{\tau(a)}{f(u,\head(a))},
	\end{align*}
	where $f(u,v)$ equals $1+|\arc_{uu}|t$ if $u=v$, $1-|\arc_{uv}||\arc_{vu}|t^2$ otherwise.

\vspace*{12pt}
\noindent
{\bf Proposition~3:} (Ishikawa, Morita and Sato \cite{ishikawa21}) 
	Let $\Delta=(V,\arc)$ be a digraph.
	The generalized weighted zeta function $Z_\Delta(t;\theta)$ is given by
	$$
		\frac{1}{\prod_{(u,v)\in\Phi_\Delta}f(u,v) \det(I-t\underline{A}^{\rm S}_\Delta+t^2 \underline{D}^{\rm S}_\Delta)}.
	$$

\vspace*{12pt}
\noindent

The matrices $A^{\rm S}_G$ and $\underbar{{\it A}}^{\rm S}_\Delta$ are called the {\it weighted adjacency matrices},
and $D^{\rm S}_G$ and $\underbar{{\it D}}^{\rm S}_\Delta$ are called the {\it weighted adjacency matrices}.

\subsection{Szegedy walk}
First we give the definition of the Szegedy walk.

\vspace*{12pt}
\noindent
{\bf Definition~2:}
For the symmetric digraph $\Delta=(V,\arc(G))$ of a simple graph $G$,
let $p : \arc(G) \to (0,1]$ be a transition probability satisfying 
$
	\sum_{a\in \arc(G), \tail(a)=v} p(a)=1
$
for each $v\in V$.
The {\it Szegedy walk} \cite{szegedy2004quantum} is the quantum walk whose transition matrix 
$U_{\rm SZ} = (u_{\rm SZ}(a,a'))_{a,a'\in \arc(G)}$ is defined by
$$
	u_{\rm SZ}(a,a') := 2\sqrt{p(a)p(\overline{a'})}\delta_{\tail(a)\head(a')}-\delta_{a' \overline{a}}.
$$

For the weight (\ref{def:weight}),
let $\tau_1(a)=\sqrt{2p(\overline{a})}$ and $\tau_2(a)=\sqrt{2p(a)}$ for each $a\in\arc(G)$.
Then, we have $\theta(a,a')=u_{\rm SZ}(a',a)$.

\vspace*{12pt}
\noindent
{\bf Remark~1:}
The {\it Grover walk} \cite{grover1966fast} is a special case of the Szegedy walk, whose transition matrix is 
$
	U_{\rm GR}= ( \frac{2}{\deg \tail(a)}\delta_{\tail(a)\head(a')} - \delta_{a' \overline{a}})_{a,a'\in \arc(G)}.
$
It is given by substituting $p=\frac{2}{\deg \tail(a)}$ for each $u_{\rm SZ}(a,a')$.

\vspace*{12pt}
\noindent


Konno and Sato \cite{konno2012relation} give the spectrum of the transition matrix of the {\it Grover walk} on $G$
by the Ihara expression of $Z_{\Delta(G)}(t;\theta)$.

\vspace*{12pt}
\noindent
{\bf Theorem~1:}\label{thm:Grover}
	For a finite simple connected graph $G$,
	\begin{align*}
		\det(\lambda I- U_{\rm GR}) &= (\lambda^2-1)^{|E|-|V|} \det((\lambda^2+1)I-2\lambda T)\\
		&= (\lambda^2-1)^{|E|-|V|} \prod_{\mu \in {\rm Spec}(T)} ((\lambda^2+1)-2\mu\lambda), \nonumber
	\end{align*}
	where $T=(T_{uv})_{u,v\in V}$ is defined as follows:
	$$
		T_{uv} =\begin{cases}
			\frac{1}{\deg(u)} & \mbox{ if } \  \{ u,v \}\in E,\\
			0 & \mbox{ otherwise}.
		\end{cases}
	$$
	Thus, we get
	$$
		{\rm Spec}(U_{\rm GR}) = \{ \pm 1 \}^{|E|-|V|} \sqcup \{ \lambda \ | \ \lambda^2 -2\mu\lambda +1=0,  \ {\rm for \ each \ }\mu\in{\rm Spec}(T) \}.
	$$

\vspace*{12pt}
\noindent

Recall that $\theta^{\rm S}(a,a') =\theta(a,a')|_{\tau_1\equiv 1} = \tau_2(a')\delta_{\head(a)\tail(a')}-\delta_{a'\in a\inv}$ for each two arcs $a,a'\in\arc$.
For each $a\in\arc(G)$,
substituting $\frac{2}{\deg \tail(a)}$ for $\tau_2(a)$ of $\theta^{\rm S}$,
we get the above theorem by the Ihara expression of $Z_{\Delta(G)}(t;\theta^{\rm S})$.
In order to consider a similar theorem for the Szegedy walk on a finite digraph $\Delta$,
it is necessary to show that $Z_{\Delta}(t;\theta)$ has the Ihara expression.

\section{Main theorem} \label{main}
\noindent

Before stating the main theorem, we will show some lemmas.

\vspace*{12pt}
\noindent
{\bf Lemma~1:}
	For a variable $t$ and a scalar $k$,
	the following hold
	\begin{align*}
		(I+t k\mathbbm{1}_n)\inv &= I - (1+tkn)\inv tk\mathbbm{1}_n, \\
		\det(I+t k\mathbbm{1}_n) &= 1+t kn.
	\end{align*}

\vspace*{12pt}
\noindent
{\bf Proof:} 
	It is easy to see from the following:
	\begin{align*}
		(I+t k\mathbbm{1}_n)\{ I - (1+tkn)\inv tk\mathbbm{1}_n \} 
		&=  I -(1+tkn)\inv tk\mathbbm{1}_n +t k\mathbbm{1}_n - (1+tkn)\inv t^2k^2 n \mathbbm{1}_n \\
		&=  I+ tk \mathbbm{1}_n - (1+tkn)\inv (1+tkn)tk\mathbbm{1}_n \\
		&=  I.
	\end{align*}
	In the same way, we can see 
	$
		\{ I - (1+tkn)\inv tk\mathbbm{1}_n \} (I+t k\mathbbm{1}_n) = I.
	$
	
	The determinant is given as follows:
	\begin{align*}
		\begin{vmatrix}
			1+t k & tk & tk & \ldots \\
			tk & 1+t k & tk & \ldots \\
			tk & tk & 1+t k & \ldots \\
			\vdots & \vdots &\vdots  & \ddots
		\end{vmatrix}
		&=
		\begin{vmatrix}
			1+t k & tk & tk & \ldots \\
			-1 & 1 & 0 & \ldots \\
			-1 & 0 & 1 & \ldots \\
			\vdots & \vdots &\vdots  & \ddots
		\end{vmatrix}\\
		&=
		\begin{vmatrix}
			1+n(tk) & tk & tk & \ldots \\
			0 & 1 & 0 & \ldots \\
			0 & 0 & 1 & \ldots \\
			\vdots & \vdots &\vdots  & \ddots
		\end{vmatrix}\\
		&= 1+tkn.
	\end{align*}
	
$\square\,$

\vspace*{12pt}
\noindent
{\bf Lemma~2:}
	Let $M_1\in\mat(k,l;\mathbb{C})$ and $M_2\in\mat(l,k;\mathbb{C})$ and 
	$
		M := \begin{bmatrix}
			O & M_1 \\
			M_2 & O
		\end{bmatrix}
	$.
	When the matrix $I_l-t^2M_2M_1$ is nonsingular,
	$(I+tM)^{-1}$ can be written by
	$$
		\begin{bmatrix}
			(I_k-t^2M_1M_2)\inv 			& -tM_1(I_l-t^2M_2M_1)\inv\\
			-t(I_l-t^2M_2M_1)\inv M_2	& (I_l-t^2M_2M_1)\inv \\
		\end{bmatrix}.
	$$
	We also have $\det(I+tM) = \det(I_l-t^2M_2M_1) = \det(I_l-t^2M_1M_2)$.

\vspace*{12pt}
\noindent
{\bf Proof:}
	The matrix $I+tM$ is decomposed as
	\begin{eqnarray}
		I+tM=
		\begin{bmatrix}
			I & O \\
			tM_2 & I
		\end{bmatrix}
		\begin{bmatrix}
			I & O \\
			O & I-t^2M_2M_1
		\end{bmatrix}
		\begin{bmatrix}
			I & tM_1 \\
			O & I
		\end{bmatrix}. \label{eq:decomp}
	\end{eqnarray}
Taking the inverse of both sides, we get
	\begin{align*}
		(I+tM)\inv &=
		\begin{bmatrix}
			I & tM_1 \\
			O & I
		\end{bmatrix}\inv
		\begin{bmatrix}
			I & O \\
			O & I-t^2M_2M_1
		\end{bmatrix}\inv
		\begin{bmatrix}
			I & O \\
			tM_2 & I
		\end{bmatrix}\inv \\
		&=
		\begin{bmatrix}
			I & -tM_1 \\
			O & I
		\end{bmatrix}
		\begin{bmatrix}
			I & O \\
			O & (I-t^2M_2M_1)\inv
		\end{bmatrix}
		\begin{bmatrix}
			I & O \\
			-tM_2 & I
		\end{bmatrix}\\
		&=
		\begin{bmatrix}
			I+tM_1(I-t^2M_2M_1)\inv tM_2 & -tM_1(I-t^2M_2M_1)\inv \\
			-(I-t^2M_2M_1)\inv tM_2 & (I-t^2M_2M_1)\inv
		\end{bmatrix}.
	\end{align*}
	By the same way as in Lemma 1,
	it can be shown that $I+tM_1(I-t^2M_2M_1)\inv tM_2 = (I_k-t^2M_1M_2)\inv$ holds.

	The matrix $I+tM$ is also decomposed as
	\begin{eqnarray}
		I+tM=
		\begin{bmatrix}
			I & tM_1 \\
			O & I
		\end{bmatrix}
		\begin{bmatrix}
			I-t^2M_1M_2 & O \\
			O & I
		\end{bmatrix}
		\begin{bmatrix}
			I & O \\
			tM_2 & I
		\end{bmatrix}. \label{eq:decomp2}
	\end{eqnarray}
	Taking the determinants of both sides of Equation (\ref{eq:decomp}) and (\ref{eq:decomp2}),
	we get 
	$$
		\det(I+tM) = \det(I-t^2M_1M_2) = \det(I-t^2M_2M_1).
	$$
$\square\,$

\vspace*{12pt}
\noindent
{\bf Remark~2:}
	For conformable matrices $A,B,C$, and $D$,
	the identity
	$
		(A+BCD)\inv = A\inv -A\inv B(C\inv +DA\inv B)\inv DA\inv
	$
	is called {\it the Woodbury matrix identity}.

\vspace*{12pt}
\noindent


\subsection{The Ihara expression of the graph zeta function on a finite digraph}
\noindent


We fix a total order $<$ on $V$,
and if one writes $\arc(u,v)$ then the condition $u<v$ is always assumed.
Let
$J=(j_{aa'})_{a,a'\in\arc}$,
$K=(k_{av})_{a\in\arc ,v\in V}$ and
$L=(l_{ua'})_{u\in V, a'\in\arc}$
be matrices defined by 
$j_{aa'}=\delta_{a'\in a\inv}$,
$k_{av}=\tau_1(a)\delta_{\head(a)v}$ and
$l_{ua'}=\tau_2(a')\delta_{u\tail(a')}$
For each $(u,v)\in\Phi_\Delta$,
let 
$J(u,v) := (j_{aa'})_{a,a'\in\arc(u,v)}$,
$K(u,v) := (k_{aw})_{a\in\arc(u,v),w\in V}$ and  
$L(u,v) := (l_{wa'})_{w\in V,a'\in\arc(u,v)}$.
Note that one can choose a total order on $V$ which makes $J$ a diagonal matrix.
We fix such a total order on $V$.
Let $T$ denote the block diagonal matrix $I+tJ$, 
and the diagonal blocks are given by $T(u,v) := I+tJ(u,v)$.
We denote $\det(T(u,v))$ by $f(u,v)$.
Then, we can see that $\det(T) = \prod_{(u,v)\in\Phi_\Delta}  f(u,v)$ holds.

\vspace*{12pt}
\noindent
{\bf Lemma~3:}
	For $(u,v)\in \Phi_\Delta$,
	\begin{align*}
		f(u,v) 
		= \begin{cases}
			1-|\arc_{uv}| |\arc_{vu}| t^2 & \mbox{ if }\  u\neq v, \\
			1+|\arc_{uu}|t &  \mbox{ if }\  u=v.
		\end{cases}
	\end{align*}

\vspace*{12pt}
\noindent
{\bf Proof:} 
	For $(u,v)\in\Phi_\Delta$ satisfying $u\neq v$, we assume that $|\arc_{uv}| =k$ and $|\arc_{vu}| =l$.
	The matrix $J(u,v)$ is given by
	$
		\begin{bmatrix}
			O & \mathbbm{1}_{k\times l} \\
			\mathbbm{1}_{l\times k} & O
		\end{bmatrix}.
	$
	Thus, we obtain
	$
		T(u,v)=
		\begin{bmatrix}
			I_k & t\mathbbm{1}_{k\times l} \\
			t\mathbbm{1}_{l\times k} & I_l
		\end{bmatrix}.
	$
	By Lemma 2,
	we have 
	$
		\det(T(u,v)) = \det(I-t^2\mathbbm{1}_{k\times l}\mathbbm{1}_{l\times k})
	$.
	Since 
	$
		\mathbbm{1}_{k\times l}\mathbbm{1}_{l \times k} = l\mathbbm{1}_{k\times k}
	$ 
	holds, 
	and the identity
	$
		\det(I-t^2 l\mathbbm{1}_{k\times k}) = 1-klt^2
	$
	follows from Lemma 1.
	
	For $(u,u)\in\Phi_\Delta$, we assume that $|\arc_{uu}| =n$.
	Since the matrix $J(u,u)$ is equal to $\mathbbm{1}_n$,
	$T(u,u)=I+t\mathbbm{1}_n$ holds.
	From Lemma 1,
	we obtain $\det(I+t\mathbbm{1}_n)=1+nt$.
$\square\,$

We consider the following matrices
\begin{align*}
	A^\theta_\Delta &= \sum_{(u,v)\in \Phi_\Delta} L(u,v) K(u,v), \\
	\underline{D}^\theta_\Delta &= \sum_{(u,v)\in \Phi_\Delta}  \frac{L(u,v) J(u,v) K(u,v)}{f(u,v)}, \\
	\underline{X}^\theta_\Delta &= \sum_{(u,v)\in \Phi_\Delta} \overline{\delta}_{uv} \frac{L(u,v) J(u,v)^2 K(u,v)}{ f(u,v)}.
\end{align*}

\vspace*{12pt}
\noindent
{\bf Theorem~2:}\label{thm:main1}
	The following identity holds:
	$$
		\det(I-tM_\theta)
		= \prod_{u,v \in V}f(u,v)  \det(I-tA^\theta_\Delta +t^2 \underline{D}^\theta_{\Delta} -t^3\underline{X}^\theta_{\Delta}).
	$$

\vspace*{12pt}
\noindent
{\bf Proof:} 
Let 
$
	H := (\tau(a,a')\delta_{\head(a)\tail(a')})_{a,a'\in\arc}
$.
Then we have $M_\theta = H-J$.
Since $\tau(a,a')\delta_{\head(a)\tail(a')}=\sum_{v\in V}(\tau_1(a)\delta_{\head(a)v})(\tau_2(a')\delta_{v\tail(a')})$,
we have $H=KL$.
Thus we obtain
$$
	\det(I-tM_\theta) =\det(I-t(KL-J)) = \det(T-tKL).
$$
For two conformable matrices $X$ and $Y$, 
it is known that $\det (I-XY) =\det (I-YX)$ holds.
Hence we have 
\begin{align*}
	\det(T-tKL) &= \det(T)\det(I-tT\inv KL) \\
	&=\left( \prod_{(u,v)\in \Phi_\Delta}f(u,v)\right) \det (I-tLT\inv K).
\end{align*}
For the direct sum decomposition $T=\bigoplus_{(u,v)\in \Phi_\Delta}T(u,v)$,
we arrange the submatrices $J(u,v)$ and $K(u,v)$ of $J$ and $K$ in order of submatrices $T(u,v)$ of $T$.
Then, we can give the matrix $LT\inv K$ by a sum
$$
	\sum_{(u,v)\in \Phi_\Delta} L(u,v)T(u,v)\inv K(u,v).
$$
The matrix $T(u,v)$ is different for $u=v$ and for $u\neq v$.
Thus we consider these two cases separately.

For $(u,u)\in\Phi_\Delta$,
we assume $|\arc(u,u)|=n$. 
Since any two arcs $a,a'\in\arc(u,u)$ are inverses of each other,
$J(u,u)=\mathbbm{1}_n$ and $T(u,u)=I+ t\mathbbm{1}_n$ holds.
By Lemma 1 and Lemma 3,
we have
\begin{align*}
	L(u,u)T(u,u)\inv K(u,u) &= L(u,u) \left( I-\frac{t}{1+nt}\mathbbm{1}_n \right) K(u,u) \\
	&= L(u,u) K(u,u) - t\frac{L(u,u)J(u,u)K(u,u)}{f(u,u)}.
\end{align*}

For $(u,v)\in\Phi_\Delta$ satisfying $u\neq v$,
we arrange the elements of $\arc(u,v)$ by an order $<$
such that $a<a'$ for $a\in\arc_{uv}$ and $a'\in\arc_{vu}$.
We assume that $|\arc_{uv}|=k$ and $|\arc_{vu}|=l$.
Then, $T(u,v)=
	\begin{bmatrix}
		I & t\mathbbm{1}_{k\times l} \\
		t\mathbbm{1}_{l\times k} & I
	\end{bmatrix}.
$
From Lemma 2, 
$(1,1)$-block and $(2,2)$-block of $T(u,v)\inv$  are 
$
	(I_k-t^2\mathbbm{1}_{k\times l} \mathbbm{1}_{l\times k} )\inv = (I_k-t^2 l \mathbbm{1}_{k\times k})\inv
$ 
and 
$
	(I_l-t^2\mathbbm{1}_{l\times k} \mathbbm{1}_{k\times l} )\inv = (I_l-t^2 k \mathbbm{1}_{l\times l})\inv
$.
Thus, $T(u,v)\inv$ is given by
$$
	\begin{bmatrix}
		(I_k-t^2 l \mathbbm{1}_{k\times k} )\inv 		& -t\mathbbm{1}_{k\times l}(I_l-t^2 k \mathbbm{1}_{l\times l})\inv\\
		-t (I_l-t^2 k \mathbbm{1}_{l\times l} )\inv \mathbbm{1}_{l\times k}	 & (I_l-t^2 k \mathbbm{1}_{l\times l})\inv
	\end{bmatrix}.
$$
By Lemma 1 and Lemma 3, we have
\begin{align*} 
	T(u,v)\inv
	&= \begin{bmatrix}
		I_k +\frac{t^2}{1-kl t^2} l\mathbbm{1}_k  & -t\mathbbm{1}_{k\times l} - \frac{kl t^3}{1-kl t^2} \mathbbm{1}_{k\times l} \\
		-t\mathbbm{1}_{l\times k} -\frac{kl t^3}{1-kl t^2}\mathbbm{1}_{l\times k} & I_l+ \frac{t^2}{1-kl t^2} k\mathbbm{1}_l
	\end{bmatrix} \\
	&=I_{k+l}
	- \frac{t}{1-kl t^2}
	\begin{bmatrix}
		O& (1-kl t^2)\mathbbm{1}_{k\times l} + kl t^2\mathbbm{1}_{k\times l} \\
		(1-kl t^2)\mathbbm{1}_{l\times k} +kl t^2\mathbbm{1}_{l\times k} &O
	\end{bmatrix}\\
	& +
	\frac{t^2}{1-kl t^2}
	\begin{bmatrix}
		 l\mathbbm{1}_k  &O\\
		O &  k\mathbbm{1}_l
	\end{bmatrix}\\
	&=I_{k+l}
	- \frac{t}{1-kl t^2}
	\begin{bmatrix}
		O&\mathbbm{1}_{k\times l}  \\
		\mathbbm{1}_{l\times k} &O
	\end{bmatrix}
	+
	\frac{t^2}{1-kl t^2}
	\begin{bmatrix}
		O&\mathbbm{1}_{k\times l}  \\
		\mathbbm{1}_{l\times k} &O
	\end{bmatrix}^2\\
	&= I_{k+l} - t \frac{J(u,v)}{f(u,v)}  + t^2 \frac{J(u,v)^2}{f(u,v)}.
\end{align*}
Hence $L(u,v) T(u,v)\inv K(u,v)$ equals 
$$
	 L(u,v) K(u,v)   -t \frac{L(u,v) J(u,v) K(u,v)}{f(u,v)} + t^2 \frac{L(u,v) J(u,v)^2 K(u,v)}{f(u,v)}.
$$

Therefore, $LT\inv K$ is given by
\begin{eqnarray*}
	\sum_{(u,v)\in\Phi_\Delta} 
	\left\{ L(u,v) K(u,v)  - t \frac{L(u,v) J(u,v) K(u,v)}{f(u,v)}  + t^2 \overline{\delta}_{uv} \frac{L(u,v) J(u,v)^2 K(u,v)}{f(u,v)} \right\}.
\end{eqnarray*}
It follows from the definition that $LT\inv K = A^\theta_\Delta -t\underline{D}^\theta_\Delta + t^2 \underline{X}^\theta_\Delta$,
which completes the proof. $\square\,$

\vspace*{12pt}
\noindent
{\bf Remark~3:}
	The entries of $A_\Delta^\theta=(a_{uv})_{u,v\in V}, \underline{D}_\Delta^\theta=(\underline{d}_{uv})_{u,v\in V}$ 
	and $\underline{X}_\Delta^\theta=(\underline{x}_{uv})_{u,v\in V}$ are given by
	\begin{align*}
		a_{uv}&= \sum_{a\in\arc_{uv}} \tau(a,a),\\
		\underline{d}_{uv}&= \delta_{uv} \sum_{w\in V} \sum_{a\in\arc_{uw}, a'\in \arc_{wu}} \frac{\tau(a', a)}{f(u,w)},\\
		\underline{x}_{uv}&= 
		 \overline{\delta}_{uv} \frac{|\arc_{vu}|}{f(u,v)} \sum_{a,a'\in\arc_{uv}}\tau(a, a').
	\end{align*}


\subsection{The Ihara expression of the graph zeta function on a finite graph}
\noindent


In order to give the characteristic polynomial of the Szegedy transition matrix by our graph zeta function,
we derive the Ihara expression of $Z_{\Delta(G)}(t;\theta)$ for a graph $G$.
Accordingly, we adopt the ``usual" definition of inverse for symmetric digraphs.

Let $G$ be a graph,
and $\Delta(G)=(V,\arc(G))$ the symmetric digraph for $G$.
Let $\arc(e)$ denote the set of two arcs induced by $e\in E$.
The two arcs in each $\arc(e)$ are inverses of each other.
Let $J(e):=(\delta_{a'\in a\inv})_{a,a'\in\arc(e)}, K(e):=(\tau_1(a)\delta_{\head(a)v})_{a\in\arc(e), v\in V}$ 
and $L(e):=(\tau_2(a')\delta_{v\tail(a')})_{v\in V,a'\in\arc(e)}$.
Note that $J(e)$ is the square matrix $\begin{bmatrix} 0 & 1 \\ 1 & 0 \end{bmatrix}$ for any $e\in E$,
and $J=\oplus_{e\in E}J(e)$ holds.
Let $T(e):=I+tJ(e)$ and $T=\oplus_{e\in E}T(e)$,
then we can see that $\det(T(e))=1-t^2$ and $\det(T)=(1-t^2)^{|E|}$.
We consider the following matrices:
\begin{eqnarray*}
	A_G^\theta := \sum_{e\in E} L(e)K(e),  \ \
	D_G^\theta := \sum_{e\in E} L(e)J(e)K(e).
\end{eqnarray*}

\vspace*{12pt}
\noindent
{\bf Theorem~3:}
	The following identity holds:
\begin{align} \label{idnt:Ihara}
	Z_{\Delta(G)}(t;\theta)\inv 
	=  (1-t^2)^{|E|-|V|}\det(I-tA_G^\theta +t^2(D_G^\theta -I)). 
\end{align}

\vspace*{12pt}
\noindent
{\bf Proof:} 

	As in Theorem 2, $\det(I-tM_\theta)$ is given by 
	$$
		\det(T-tKL) = \det(T)\det(I-tLT\inv K) = (1-t^2)^{|E|} \det(I-tLT\inv K).
	$$
	Since $T=\bigoplus_{e\in E}T(e)$, we have
	$
		LT\inv K = \sum_{e\in E} L(e)T(e)\inv K(e)
	$.
	For each $e\in E$,
	$
		T(e)\inv = (1-t^2)\inv \begin{bmatrix} 1 & -t \\ -t & 1 \end{bmatrix} 
		=  (1-t^2)\inv (I-tJ(e)),
	$
	and we get
	$$
		L(e)T(e)\inv K(e) = (1-t^2)\inv \{ L(e) K(e) -L(e)J(e) K(e) \}.
	$$
	Since $LT\inv K=\sum_{e\in E} (1-t^2)\inv \{ L(e) K(e) -L(e)J(e) K(e) \}$ holds, we see
	\begin{align*}
		&(1-t^2)^{|E|}\det(I-tLT\inv K) \\ 
		&= (1-t^2)^{|E|-|V|}\det((1-t^2)I-t \sum_{e\in E} ( L(e) K(e) -L(e)J(e) K(e) )) \\
		&= (1-t^2)^{|E|-|V|}\det(I-t A_G^\theta +t^2(D_G^\theta -I)).
	\end{align*}
 $\square\,$

\vspace*{12pt}
\noindent
{\bf Remark~4:}
	The entries of $A_G^\theta=(a_{uv})_{u,v\in V}$ and $D_G^\theta=(d_{uv})_{u,v\in V}$ are given by
	\begin{align*}
		a_{uv} = \sum_{a\in\arc_{uv}} \tau(a,a), \ \ \ 
		d_{uv}= \delta_{uv} \sum_{a\in\arc_{u*}} \tau(\overline{a}, a).
	\end{align*}
\vspace*{12pt}
\noindent

As was mentioned in the introduction,
we obtain a generalization of the Konno-Sato's theorem \cite{konno2012relation}.

\vspace*{12pt}
\noindent
{\bf Corollary~1:}
	For a graph $G=(V,E)$ without loops,
	the characteristic polynomial of the transition matrix of the Szegedy walk on $G$
	is given as follows:
	\begin{align*}
		\det(\lambda I-U_{\rm SZ}) &= (\lambda^2-1)^{|E|-|V|} \left( (\lambda^2+1) I - 2\lambda T \right) \\
		&=	 (\lambda^2-1)^{|E|-|V|} \prod_{\mu\in {\rm Spec}(T)}\left( (\lambda^2+1) - 2\mu \lambda \right), 
	\end{align*}
	where
	$
		T=\left( \sum_{a\in\arc_{uv}} \sqrt{p(a) p(\overline{a})}\right)_{u,v\in V}.
	$ 
	
\vspace*{12pt}
\noindent

\vspace*{12pt}
\noindent
{\bf Proof:}
	Let $\tau_1(a)=\sqrt{2p(\overline{a})}$ and $\tau_2(a)=\sqrt{2p(a)}$ for $a\in\arc(G)$.
	Then $\theta(a,a')=u_{\rm SZ}(a',a)$ and $M_\theta= ^tU_{\rm SZ}$ hold.
	The matrices $A_G^\theta$ and $D_G^\theta$ of Equation (\ref{idnt:Ihara}) are given by
	\begin{align*}
		A_G^\theta = \left( \sum_{a\in\arc_{uv}} 2\sqrt{p(a)p(\overline{a})} \right)_{u,v\in V} = 2T, \ \ 
		D_G^\theta = \left(\delta_{uv} \sum_{a\in\arc_{u*}} 2p(a) \right)_{u,v\in V} = 2I_{|V|}.
	\end{align*}
	Thus, we obtain
	\begin{align*}
		\det(\lambda I-U_{\rm SZ}) &= \det(\lambda I- ^tU_{\rm SZ}) \\
		&= \det(\lambda I-M_\theta) \\
		&= (\lambda^2-1)^{|E|-|V|}\det(\lambda^2I -\lambda A_G^\theta +(D_G^\theta -I)) \\
		&= (\lambda^2-1)^{|E|-|V|}\det(\lambda^2I -2\lambda T + I).
	\end{align*}

 $\square\,$

\section{Example} \label{sec:exam}
\noindent

\begin{figure}[h]
\center
\includegraphics[width=6cm]{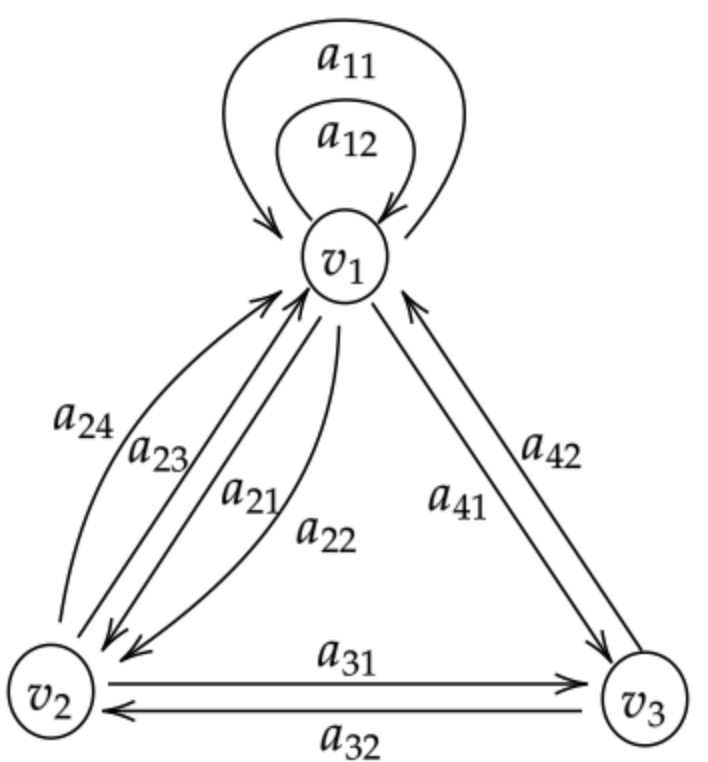}
\label{fig:digraph}
\caption{digraph $\Delta=(V,\arc)$}
\end{figure}

Let $\Delta=(V,\arc)$ be a digraph with $V=\{ v_1,v_2,v_3 \}$ and 
$\arc=\{ 
	a_1=( v_1,v_1) < a_2=( v_1,v_1)  < 
	a_3=( v_1,v_2) < a_4=( v_1,v_2) <  a_5=( v_2,v_1) <  a_6=( v_2,v_1)  < 
	a_7=(v_2,v_3) < a_8=(v_3,v_2) < 
	a_9=( v_1,v_3) < a_{10}=( v_3,v_1) 
\}$ (see Figure \ref{fig:digraph}).
The matrices $J,K,$ and $L$ are as follows:
\begin{eqnarray*}\small
	J=
	\begin{bmatrix}
		1&1&&&&&&&& \\
		1&1&&&&&&&& \\
		&&0&0&1&1&&&& \\
		&&0&0&1&1&&&& \\
		&&1&1&0&0&&&& \\
		&&1&1&0&0&&&& \\
		&&&&&&0&1&& \\
		&&&&&&1&0&& \\
		&&&&&&&&0&1 \\
		&&&&&&&&1&0 
	\end{bmatrix}, \ \
	K=
	\begin{bmatrix}
		\tau_1(a_1)&& \\
		\tau_1(a_2)&& \\
		&\tau_1(a_3)& \\
		&\tau_1(a_4)& \\
		\tau_1(a_5)&& \\
		\tau_1(a_6)&& \\
		&&\tau_1(a_7) \\
		&\tau_1(a_8)& \\
		&&\tau_1(a_9)  \\
		\tau_1(a_{10}) &&
	\end{bmatrix}, \\
	L=
	\begin{bmatrix}
		\!\tau_2(a_1)\!\!&\!\!\tau_2(a_2)\!\!&\!\!\tau_2(a_3)\!\!&\!\!\tau_2(a_4)\!\!&&&&&\!\!\tau_2(a_9)\!\!&\!\!   \\
		&&&&\!\!\tau_2(a_5)\!\!&\!\!\tau_2(a_6)\!\!&\!\!\tau_2(a_7)\!\!&&&\!\!   \\
		&&&&&&&\!\!\tau_2(a_8)\!\!&&\!\!\tau_2(a_{10})\! 
	\end{bmatrix}.
\end{eqnarray*}
Then, $A_\Delta^\theta, \underline{D}_\Delta^\theta$ and $\underline{X}_\Delta^\theta$ are
\begin{align*}
	&A_\Delta^\theta = 
	\begin{bmatrix}
		\tau(a_1,a_1)+\tau(a_2,a_2) & \tau(a_3,a_3)+\tau(a_4,a_4)&\tau(a_9,a_9) \\
		\tau(a_5,a_5)+\tau(a_6,a_6)&0&\tau(a_7,a_7) \\
		\tau(a_{10},a_{10})&\tau(a_8,a_8)&0
	\end{bmatrix}, \\
	&\underline{D}_\Delta^\theta = \begin{bmatrix}
		d_{v_1,v_1}+d_{v_1,v_2}+d_{v_1,v_3} &0&0 \\
		0&d_{v_2,v_1}+d_{v_2,v_3}&0 \\
		0&0&d_{v_3,v_1}+d_{v_3,v_2}
	\end{bmatrix}, \\
	&\underline{X}_\Delta^\theta = 
	\begin{bmatrix}
		0 & x_{v_1,v_2}& x_{v_1,v_3}\\
		x_{v_2,v_1}& 0 & x_{v_2,v_3}\\
		x_{v_3,v_1} & x_{v_3,v_2} & 0
	\end{bmatrix},
\end{align*}
where 
\begin{align*}
	&d_{v_1,v_1} = \frac{\tau(a_1,a_1)+\tau(a_1,a_2)+\tau(a_2,a_1)+\tau(a_2,a_2)}{1+2t}, \\
	&d_{v_1,v_2} = \frac{\tau(a_5,a_3)+\tau(a_6,a_3)+\tau(a_5,a_4)+\tau(a_6,a_4)}{1-4t^2}, \qquad
	d_{v_1,v_3} = \frac{\tau(a_{10},a_9)}{1-t^2},\\
	&d_{v_2,v_1} = \frac{\tau(a_3,a_5)+\tau(a_4,a_5)+\tau(a_3,a_6)+\tau(a_4,a_6)}{1-4t^2},\qquad
	d_{v_2,v_3} =  \frac{\tau(a_8,a_7)}{1-t^2} , \\
	&d_{v_3,v_1} =  \frac{\tau(a_9,a_{10})}{1-t^2}, \qquad 
	d_{v_3,v_2} =  \frac{\tau(a_7,a_8)}{1-t^2}, \\
	&x_{v_1,v_2} = \frac{2(\tau(a_3,a_3)+\tau(a_3,a_4)+\tau(a_4,a_3)+\tau(a_4,a_4))}{1-4t^2},\qquad
	x_{v_1,v_3} = \frac{\tau(a_9,a_9)}{1-t^2},\\
	&x_{v_2,v_1} = \frac{2(\tau(a_5,a_5)+\tau(a_5,a_6)+\tau(a_6,a_5)+\tau(a_6,a_6))}{1-4t^2},\qquad
	x_{v_2,v_3} = \frac{\tau(a_7,a_7)}{1-t^2}, \\
	&x_{v_3,v_1} = \frac{\tau(a_{10},a_{10})}{1-t^2}, \qquad
	x_{v_3,v_2} = \frac{\tau(a_8,a_8)}{1-t^2}.
\end{align*}

\begin{figure}[h]
\center
\includegraphics[width=6cm]{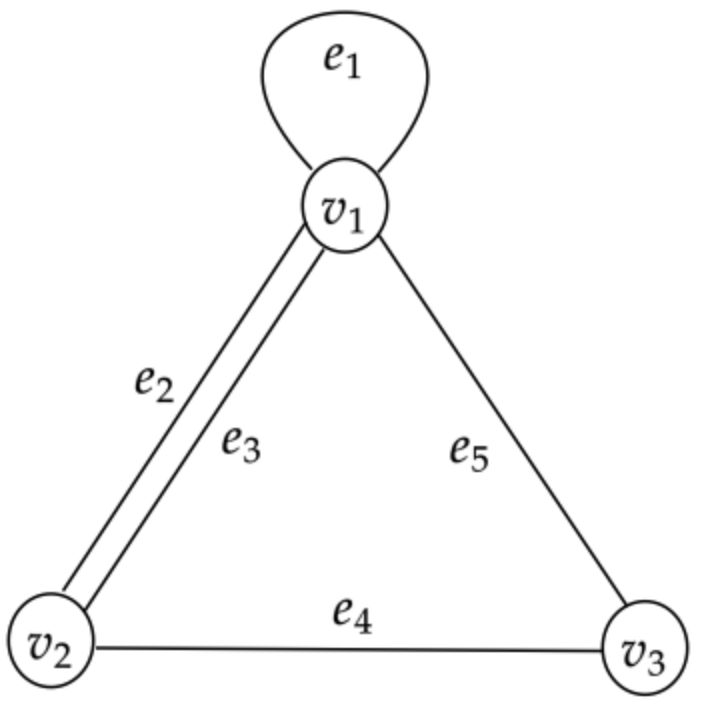}
\label{fig:graph}
 \caption{graph $G=(V,E)$}

\end{figure}

Let $G=(V,E)$ be a graph with $V=\{ v_1,v_2,v_3 \}$ and 
$E=\{ e_1=\{ v_1,v_1\}, e_2=\{ v_1,v_2\}, e_3=\{ v_1,v_2\}, e_4=\{ v_2,v_3\}, e_5=\{ v_1,v_3\} \}$.
Assign 
$e_1$ to $\{ a_1,a_2 \}$, 
$e_2$ to $\{ a_3,a_5 \}$, 
$e_3$ to $\{ a_4,a_6 \}$, 
$e_4$ to $\{ a_7,a_8 \}$ and
$e_5$ to $\{ a_9,a_{10} \}$,
where the directions of $a_1,\ldots, a_{10}$ are the same as in $\Delta$.
Thus, the symmetric digraph $\Delta(G)$ is the same as $\Delta$.
Note that the two arcs corresponding to the same edge are inverses of each other, 
and they do not have any other inverses.
Let $\Delta(G)=(V,\arc(G))$ be the symmetric digraph of $G$
with
$\arc(G)=\{ a_1< a_2 < a_3<a_5<a_4<a_6<a_7<a_8 <a_9<a_{10} \}$.
The matrices $J,K,$ and $L$ are as follows:
\begin{eqnarray*}\small
	J=
	\begin{bmatrix}
		0&1&&&&&&&& \\
		1&0&&&&&&&& \\
		&&0&1&&&&&& \\
		&&1&0&&&&&& \\
		&&&&0&1&&&& \\
		&&&&1&0&&&& \\
		&&&&&&0&1&& \\
		&&&&&&1&0&& \\
		&&&&&&&&0&1 \\
		&&&&&&&&1&0 
	\end{bmatrix}, \ \
	K=
	\begin{bmatrix}
		\tau_1(a_1)&& \\
		\tau_1(a_2)&& \\
		&\tau_1(a_3)& \\
		\tau_1(a_5)&& \\
		&\tau_1(a_4)& \\
		\tau_1(a_6)&& \\
		&&\tau_1(a_7) \\
		&\tau_1(a_8)& \\
		&&\tau_1(a_9)  \\
		\tau_1(a_{10}) &&
	\end{bmatrix}, \\
	L=
	\begin{bmatrix}
		\!\tau_2(a_1)\!\!&\!\!\tau_2(a_2)\!\!&\!\!\tau_2(a_3)\!\!&&\!\!\tau_2(a_4)\!\!&&&&\!\!\tau_2(a_9)\!\!&\!\!   \\
		&&&\!\!\tau_2(a_5)\!\!&&\!\!\tau_2(a_6)\!\!&\!\!\tau_2(a_7)\!\!&&&\!\!   \\
		&&&&&&&\!\!\tau_2(a_8)\!\!&&\!\!\tau_2(a_{10})\! 
	\end{bmatrix}.
\end{eqnarray*}
Then, $A_\Delta^\theta$ and $D_\Delta^\theta$ are 
\begin{eqnarray*} 
	&&A_\Delta^\theta = 
	\begin{bmatrix}
		\tau(a_1,a_1)+\tau(a_2,a_2) & \tau(a_3,a_3)+\tau(a_4,a_4)&\tau(a_9,a_9) \\
		\tau(a_5,a_5)+\tau(a_6,a_6)&0&\tau(a_7,a_7) \\
		\tau(a_{10},a_{10})&\tau(a_8,a_8)&0
	\end{bmatrix}, \\
	&&D_\Delta^\theta = \begin{bmatrix}
		d_{v_1,v_1}+d_{v_1,v_2}+d_{v_1,v_3} &0&0 \\
		0&d_{v_2,v_1}+d_{v_2,v_3}&0 \\
		0&0&d_{v_3,v_1}+d_{v_3,v_2}
	\end{bmatrix}, \\
\end{eqnarray*}
where 
\begin{align*}
	&d_{v_1,v_1} =  \tau(a_1,a_2)+\tau(a_2,a_1), \quad 
	d_{v_1,v_2} = \tau(a_5,a_3)+\tau(a_6,a_4),\quad 
	d_{v_1,v_3} = \tau(a_{10},a_9), \\
	&d_{v_2,v_1} = \tau(a_3,a_5)+\tau(a_4,a_6),\quad
	d_{v_2,v_3} = \tau(a_8,a_7), \\
	&d_{v_3,v_1} = \tau(a_9,a_{10}), \quad
	d_{v_3,v_2} =  \tau(a_7,a_8).
\end{align*}

\section*{Acknowledgements}
\noindent
A. Ishikawa is partially supported by Grant-in-Aid for JSPS Fellows (Grant No. 20J20590).



\noindent
\begin{thebibliography}{000}
\bibitem{bartholdi1999counting}
L. Bartholdi (1999), 
{\it Counting paths in graphs},
Enseign. Math, Vol. 45,
{\rm pp. 83-131.}

\bibitem{grover1966fast}
Lov K. Grover (1996), 
{\it A fast quantum mechanical algorithm for database search},
Proceedings of the Twenty-Eighth Annual ACM Symposium on Theory of Computing,
{\rm pp. 212-219.}

\bibitem{IDE2021227}
Y. Ide, A. Ishikawa, H. Morita, I. Sato and E. Segawa (2021), 
{\it The Ihara expression for the generalized weighted zeta function of a finite simple graph},
Linear Algebra and its Applications, Vol. 627,
{\rm pp. 227-241.}

\bibitem{ishikawa21}
A. Ishikawa, H. Morita and I. Sato,
{\it The Ihara expression for generalized weighted zeta functions of Bartholdi type on finite digraphs},
in preparation.

\bibitem{ihara1966discrete}
Y. Ihara (1966),
{\it On discrete subgroups of the two by two projective linear group over p-adic fields},
Journal of the Mathematical Society of Japan, 
Vol. 18, 
{\rm pp.219-235.}

\bibitem{konno2012relation}
N. Konno and I. Sato (2012),
{\it On the relation between quantum walks and zeta functions},
Quantum Information Processing, 
Vol. 11,
{\rm pp.341-349.}

\bibitem{konno2017}
N. Konno, I. Sato and E. Segawa (2017),
{\it The spectra of the unitary marix of a 2-tessellable staggered quantum walk on a graph},
Yokohama Math. J., 
Vol. 62,
{\rm pp.52-87.}

\bibitem{matsuoka2011theoretical}
L. Matsuoka, A. Ichihara, M. Hashimoto, and K. Yokoyama (2011),
{\it Theoretical study for laser isotope separation of heavy-element molecules in a thermal distribution},
In Proceedings of the International Conference Toward and Over the Fukushima Daiichi Accident (GLOBAL 2011), 
392063.

\bibitem{morita2020ruelle}
H. Morita (2020),
{\it Ruelle zeta functions for finite digraphs},
Linear Algebra and its Applications, 
Vol. 603,
{\rm pp.329-358.}

\bibitem{Orrell_2020}
D. Orrell (2020),
{\it A quantum walk model of financial options},
SSRN Electronic Journal.

\bibitem{Portugal_2018}
R. Portugal (2018),
{\it Quantum Walks and Search Algorithms},
Springer International Publishing, 2nd edition.

\bibitem{sato2007new}
 I. Sato (2007),
 {\it A new Bartholdi zeta function of a graph},
 International Journal of Algebra, 
 Vol. 1,
 {\rm pp.269-281.}
 
\bibitem{szegedy2004quantum}
M. Szegedy (2004),
{\it Quantum speed-up of markov chain based algorithms},
In 45th Annual IEEE symposium on foundations of computer science, 
{\rm pp.32-41.}

\end{thebibliography}
\end{document}